\documentclass{article}
\usepackage{maa-monthly}

\theoremstyle{theorem}
\newtheorem{theorem}{Theorem}
\newtheorem{prop}[theorem]{Proposition}
 
\theoremstyle{definition}
\newtheorem{definition}{Definition}
\newtheorem*{remark}{Remark}

\begin{document}

\title{Chords of an Ellipse, Lucas Polynomials, \newline and Cubic Equations}
\markright{Chords, Lucas polynomials, and cubics}
\author{Ben Blum-Smith and Japheth Wood}

\maketitle

\begin{abstract}
A beautiful theorem of Thomas Price links the Fibonacci numbers and the Lucas polynomials to the plane geometry of an ellipse, generalizing a classic problem about circles. We give a brief history of the circle problem, an account of Price's ellipse proof, and a reorganized proof, with some new ideas, designed to situate the result within a dense web of connections to classical mathematics. It is inspired by Cardano's solution of the cubic equation and Newton's theorem on power sums, and yields an interpretation of generalized Lucas polynomials in terms of the theory of symmetric polynomials. We also develop additional connections that surface along the way; e.g., we give a parallel interpretation of generalized Fibonacci polynomials, and we show that Cardano's method can be used write down the roots of the Lucas polynomials.
\end{abstract}

\section{Introduction.}

A classic problem instructs the reader to mark off $n$ equally spaced points on a unit circle, draw the chords connecting one of the points to the remaining $n-1$ others, and find the product of the lengths.

\begin{figure}[htbp]\label{fig:circlechords}
\begin{center}
\resizebox{2in}{!}{\includegraphics{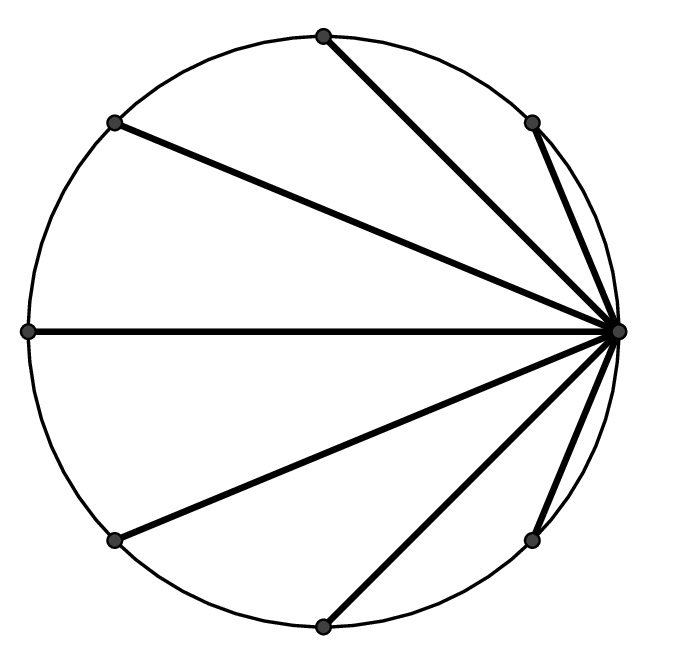}}
\caption{Chords of a circle.}
\label{diagonals}
\end{center}
\end{figure}

Although most or all of the lengths are irrational, the product is exactly $n$. An elegant solution \cite{galovich, price1, silva} involves the $n$th roots of unity $1, \zeta, \zeta^2, \ldots, \zeta^{n-1}$, which are equally spaced around the unit circle in the complex plane. Connecting $1$ to each of the $n - 1$ others and multiplying the lengths of the segments gives the product
$$\Pi_n = |1 - \zeta | |1 - \zeta^2| \cdots |1 - \zeta^{n-1}| = |(1 - \zeta)(1 - \zeta^2) \cdots (1 - \zeta^{n-1})|.$$

\noindent This is the absolute value of the polynomial $$\Pi_n(z) = (z - \zeta)(z - \zeta^2) \cdots (z - \zeta^{n-1})$$ evaluated at $z=1$.

The $n$th roots of unity are, by definition, the roots of $z^n-1$.  Now $\Pi_n$ has as roots all $n$th roots of unity except 1; it follows that $\Pi_n(z)(z-1) = z^n - 1$, and so $\Pi_n(z)=z^{n-1}+\dots+z+1$.  Thus the answer desired is $|\Pi_n(1)| = n$.

This problem has an interesting history, of which we discuss a few highlights below in Section \ref{sec:circlehistory}.

In \cite{price1}, Thomas E. Price considered the following generalization. Scale Figure \ref{fig:circlechords} horizontally and/or vertically. Then the circle becomes an ellipse. What happens to the product of the chord lengths? An elegant explicit formula is given in \cite{price1}, which we re-derive twice below (equation \eqref{eq:finalans} and Proposition \ref{prop:price}).

\begin{figure}[htbp]
\begin{center}

\resizebox{2in}{!}{\includegraphics{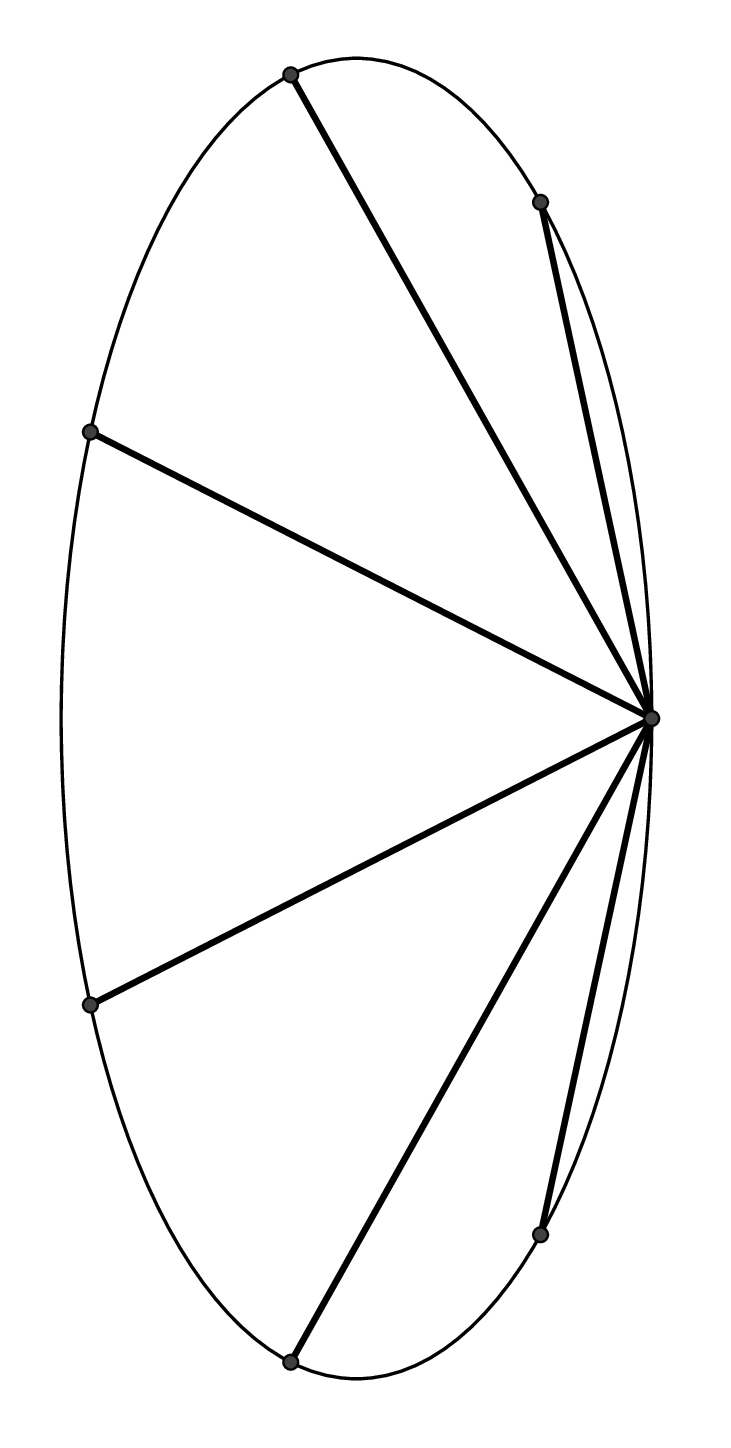}}

\caption{Stretched chords.}
\label{fig:ellipsechords}
\end{center}
\end{figure}

Showcasing the power of his results, Price in \cite{price2} considered the special case of a vertical stretch factor of $\sqrt{5}$ (see Figure \ref{fig:ellipsechords}). In this situation, the product of the chord lengths is given in Table \ref{tbl:ellipseproducts}. The answer is remarkable.

\begin{table}
\caption{Products of chords stretched by a factor of $\sqrt{5}$. The symbol $F_n$ refers to the $n$th Fibonacci number.}\label{tbl:ellipseproducts}
\begin{center}
\begin{tabular}{c|c}
$n$ & product\\ \hline
2 & $2 = 2\cdot 1$\\
3 & $6 = 3\cdot 2$\\
4 &  $12 = 4\cdot 3$\\
5 &  $25 = 5\cdot 5$\\
6 & $48 = 6\cdot 8$\\
7 & $91= 7\cdot 13$\\
\vdots & \vdots \\
$n$ & $nF_n$
\end{tabular}
\end{center}
\end{table}

The proof of Price's ellipse theorem is very clever. As is often true of very clever proofs, it is also somewhat mysterious. Our main purpose in this article is to offer a reorganized account, with some new ideas, of this beautiful and insufficiently-well-known result of Price. Our development is optimized to provide illumination and context. It aims to situate the argument within a dense web of connections to classical theory---the fundamental theorem on symmetric polynomials, Cardano's solution to the cubic in radicals, Newton's theorem on power sums, and the theories of linear recurrences and of generalized Fibonacci and Lucas polynomials all make appearances---and to use these connections to contextualize and motivate the steps. We also trace some history, explicate Price's own proof, and further develop some connections that surface.

The structure of the article is as follows. In Section \ref{sec:circlehistory}, we assemble a history (far from complete) of the theorem on circles that opened this article and provided Price with his inspiration. In this section we also distinguish five steps in the proof above, which are used both to speak precisely about the historical development, and to establish the framework in which the ellipse problem will be solved. In Section \ref{sec:outline}, we carry the five steps of the circle proof over to the ellipse, reducing the problem to the evaluation of the derivative of a certain polynomial $\Omega_n(z)$ at a certain real number $a+b$. Then we explicate Price's own proof of his beautiful formula, i.e., his evaluation of $\Omega_n'(a+b)$, and outline the rest of the article. Sections \ref{sec:cardano}--\ref{sec:solution} present our reorganized evaluation of $\Omega_n'(a+b)$. Each of these sections presents one step in the argument near the beginning---Propositions \ref{prop:thepoly}--\ref{prop:price}, respectively---and then goes on to discuss connections with known theory, motivation, and relations with Price's argument.  Section \ref{sec:furtherremarks} develops an additional connection revealed by Sections \ref{sec:cardano}--\ref{sec:solution}. A more detailed outline of the contents of Sections \ref{sec:cardano}--\ref{sec:furtherremarks} is given at the end of Section \ref{sec:outline}, after the framework has been set up.

Throughout, $\zeta$ refers to a primitive $n$th root of unity, which can be taken to be $e^{2\pi i / n}$. We prefer it to $\zeta_n$ to avoid visual clutter, but the notation conceals the dependence of $\zeta$ on $n$. Hopefully no confusion will result.

\section{History of the circle problem.}\label{sec:circlehistory}

Price's inspiration, the circle problem described in the introduction, has its own interesting print history, of which we share some highlights here. We think it is likely that the problem has been rediscovered many times, so we make no attempt to be exhaustive. As given in the introduction, the proof (henceforth, the ``standard proof") that the product of the chord lengths is $n$ (henceforth, the ``circle theorem") consists of five steps:
\begin{enumerate}
\item Interpret the $n$ equidistant points on the unit circle as the $n$th roots of unity $1,\zeta,\dots,\zeta^{n-1}$ in the complex plane $\mathbb{C}$.\label{step:demoivre}
\item Interpret the relevant chord lengths as absolute values of the complex numbers $1-\zeta^j$, $j=1,\dots,n-1$.\label{step:geom}
\item Exploit the fact that products commute with absolute values in $\mathbb{C}$ to identify the product of the chord lengths with $(1-\zeta)\cdots(1-\zeta^{n-1})$.\label{step:absval}
\item Note that this is the polynomial $\Pi_n(z) = (z-\zeta)\cdots(z-\zeta^{n-1})$, evaluated at $z=1$.\label{step:poly}
\item Invoke the identity $\left(z-\zeta\right)\cdots(z-\zeta^{n-1}) = z^{n-1} + \cdots + z + 1$ to carry out this evaluation.\label{step:identity}
\end{enumerate}
Taken individually, each of these is routine, very classical, or both. Thus a broadly-construed ``history of the circle theorem" would be a history of the complex plane, the roots of unity, and their connection to the circle. This is beyond our present scope. But we mention some striking signposts.

The first of these is a result quite close to the circle theorem that was discovered in 1716 by Roger Cotes, and enunciated (without proof) in his {\em Harmonia Mensurarum}, published posthumously in 1722 (see \cite[pp.~194--195]{stillwell}). Cotes is better known as the editor of the second edition of Newton's {\em Principia}. The following statement of Cotes's theorem is taken from \cite[p.~195]{stillwell}:

{\em If $A_0,\dots,A_{n-1}$ are equally spaced points on the unit circle with center $O$, and if $P$ is a point on $OA_0$ such that $OP=x$, then
\[
PA_0\cdot PA_1\cdot \dots \cdot PA_{n-1} = 1 - x^n.
\]}
One obtains the circle theorem by dividing through by $PA_0 = 1-x$ and then letting  $P\to A_0$ (equivalently, $x\to 1$). 

How Cotes came to this conclusion has not been preserved, but he and his contemporaries Johann Bernoulli and Abraham de Moivre were taking halting steps toward the formula
\begin{equation}\label{eq:demoivre}
(\cos \theta + i\sin\theta)^n = \cos n\theta + i\sin n\theta.
\end{equation}
This formula now bears de Moivre's name \cite[pp.~192--195]{stillwell}, although it never appeared in his writings \cite[p.~193]{stillwell}. To a modern reader, and indeed by later in the 18th century, de Moivre's formula would be seen as a consequence of the formula
\[
e^{i\theta} = \cos\theta + i\sin \theta
\]
given by Euler \cite[Chapter 8, article 138]{euler} in 1748.

John Stillwell has speculated \cite[p.~195]{stillwell} that Cotes may have used reasoning related to some variant of \eqref{eq:demoivre}. To us, this formula is the starting point of the circle theorem, for it gives the identification (step \ref{step:demoivre}) of the roots of unity with points on a unit circle: set the right side of \eqref{eq:demoivre} to $1$, obtaining $n\theta=2\pi k$ for some integer $k$, and then solve for $\theta$, yielding $\theta = 2\pi k / n$. Then consult the left side to conclude $(\cos 2\pi k / n + i\sin 2\pi k / n)^n = 1$, so that the $n$th roots of unity are the numbers $\cos 2\pi k/n + i\sin 2\pi k / n$. By the late 18th century, this was a standard maneuver---see for example \cite[article 23, p.~249]{lagrange}.

We are about to skip ahead 200 years, but we mention two major 19th century developments that are mathematically adjacent to our story. 

First, the identity of step \ref{step:identity} is the starting point for the chapter on cyclotomic (``circle-dividing") equations in Gauss's 1801 {\em Disquisitiones Arithmeticae} \cite[Chapter 7]{gauss}. In this work, Gauss proved that any root of unity $\zeta$ can be expressed in radicals, via calculations in what we now recognize as the Galois group of the polynomial $\Pi_n(z)=z^{n-1} + \dots + z + 1$. This was the first real use of Galois groups---30 years before Galois! As a corollary, he deduced that if $n$ is prime, a regular $n$-gon can be constructed with ruler and compass if and only if $n$ has the form $2^k + 1$, and therefore the 17-gon is constructible.

Second, the special case of the identity in step \ref{step:identity} with the substitution $z=1$ (as in steps \ref{step:absval}, \ref{step:poly}) figures in the series of monumental papers by Ernst Kummer \cite{kummerw} that birthed the field of algebraic number theory. Kummer investigated how factorization changes when one expands from the ring of integers $\mathbb{Z}$ to the bigger ring $\mathbb{Z}[\zeta]$ consisting of integer linear combinations of $\zeta$ and its powers. When $n=p$ is prime, substituting $z=1$ into the identity of step \ref{step:identity} yields $p=(1-\zeta)\cdots(1-\zeta^{p-1})$, which is the factorization of $p$ into primes in $\mathbb{Z}[\zeta]$ (e.g., \cite[p.~174]{kummer}). This work was important to Kummer's proof of Fermat's last theorem for so-called {\em regular primes}, which was the first work on Fermat's last theorem that proved the full statement for a large class of prime exponents at once. See \cite[Chapters 4--6]{edwards} for a modern exposition of Kummer's contribution.

We now reach the 20th century. The circle theorem was stated in a short 1954 note of W. Sichardt \cite{sichardt} in \textit{Zeitschrift f\"{u}r Angewandte Mathematik und Mechanik} (\textit{ZAMM}). This is the first instance of which we are aware in which the theorem was given in the form we have stated it, although again we make no claim to comprehensiveness.

Sichardt made the observation much earlier---in 1927---in the context of an engineering problem having to do with the construction of wells! The calculation of this particular product of chords was directly motivated by the engineering context.

Sichardt described finding the pattern in the products empirically. He did not claim credit for the proof. He stated that an unnamed mathematician who worked for Siemens provided a proof that was lost during the war. He attributed the proof given in the note to a Professor Szab\'{o} of the Technical University of Berlin. It is a slightly less clean version of the standard proof. It makes use of the fact that the length of a chord in the unit circle subtended by an angle $\theta$ is $2\sin(\theta/2)$, and thus the product of interest can be expressed as
\[
2^{n-1}\prod_{j=1}^{n-1} \sin \pi j/n.
\]
Each factor $\sin \pi j / n$ is expressed as $\frac{1}{2i}\left(\eta^j - \eta^{-j}\right)$, where $\eta$ is a $2n$th root of unity rather than an $n$th, so a certain amount of bookkeeping is needed to arrive at the expression in terms of the polynomial $\Pi_n(z)$. The comparative cleanness of the standard proof comes from the fact that taking absolute values allows one to forgo all this bookkeeping.

In 1972, Kurt Eisemann, in another short note in {\em ZAMM} \cite{eisemann}, extended Sichardt's result by considering the product of the lengths of perpendiculars from the center of the circle to the chords. Independently, Zalman Usiskin \cite{usiskin} in 1979 derived various identities involving products of sines using only the standard trigonometric identities $\sin 2\theta = 2\sin\theta\cos\theta$ and $\cos\theta = \sin (\pi/2-\theta)$, and noted the interpretation of these identities in terms of chord lengths. One of the identities proven by Usiskin yielded the case $n=45$ of the circle theorem, and he stated the general case without proof.

In 1987, Steven Galovich \cite{galovich} took up Usiskin's challenge to find a general principle explaining his various sine-product identities. Among other results, Galovich proved the circle theorem, using exactly the standard proof. But he introduced it with the words, ``Although the next theorem and proof are evidently well known, it is natural to include them in this note." \cite[p.~112]{galovich} 

By this point, the circle theorem had been posed as an exercise in textbooks, for example \cite[p.~16, exercise 12]{baknewman} and \cite[p.~69, problem 44]{kreyszig}.

In 1995, Andre Mazzoleni and Samuel Shan-Pu Shen published a note \cite{mazzolenishen} in the February issue of {\em Mathematics Magazine} giving a short proof of the circle theorem via the theory of residues of a function of a complex variable. For precedents for the result, they cited exercises in several textbooks in complex analysis. Responses from readers were published in the June and October issues pointing out other precedents. Among these readers were Usiskin, and also Eisemann, who called attention to Sichardt's contribution as well as his own. One reader fit an extremely concise version of the standard proof into a letter to the editor \cite{silva}.

In 2002, Barry Lewis gave a riff \cite{lewis} on the circle theorem by deriving formulas for power sums of the chord lengths, rather than their product.

Finally we come to Price's work in the early 2000's. In addition to \cite{price1, price2} which are the starting point of the present article, Price also produced \cite{price3}, extending Eisemann's \cite{eisemann} results to the elliptical situation. In these works, Price cited Sichardt for the circle theorem, and gave the standard proof.

\section{Price's ellipse argument; this article's goals.}\label{sec:outline}

In this section, we set the mathematical framework that will be used, explicate Price's solution to the ellipse problem, and outline the goals and structure of the rest of the article.

It is natural to attempt to replicate steps \ref{step:demoivre}--\ref{step:identity} of the standard proof (see the previous section) in the setting of the ellipse, and this is essentially what Price did in \cite{price1}. Suppose the unit circle is stretched horizontally by a factor of $a+b$ and vertically by $a-b$, where $a > |b| \geq 0$. A computation with real and imaginary parts then shows that this stretching sends the original marked points $\zeta^j,\; j=0,\dots,n-1$ in the complex plane to the points
\begin{equation}\label{eq:roots}
\zeta^ja + \zeta^{-j}b.
\end{equation}
(The original unit circle of Figure \ref{fig:circlechords} is the case $a=1$, $b=0$. In the case where the circle is vertically scaled by $\sqrt{5}$, as in Figure \ref{fig:ellipsechords}, we have $a=\phi = (1+\sqrt{5})/2$, the golden ratio, and $b=\hat\phi = (1-\sqrt{5})/2$, its algebraic conjugate.) 

Replicating steps \ref{step:demoivre}--\ref{step:poly} in the new situation is then straightforward: exactly the same interpretations and manipulations take place, but with $\zeta^j a + \zeta^{-j}b$ in the place of $\zeta^j$ everywhere, and, in particular, $a+b$ in the place of $1=\zeta^0$. Thus, in step \ref{step:poly}, the product of elliptical chord lengths is identified as the absolute value of the polynomial
\[
\Pi_n^{\text{ellipse}}(z) = \prod_{j=1}^{n-1} \left( z - (\zeta^{j}a + \zeta^{-j}b)\right),
\]
evaluated at $z=a+b$. 

The new work lies in replicating step \ref{step:identity}. One needs a way to carry out the evaluation of $\Pi_n^\text{ellipse}(a+b)$. Price's \cite{price1} elegant method proceeds as follows:

First, define
\[
\Omega_n(z)= (z-(a+b))\Pi_n^{\text{ellipse}}(z),
\]
the unique monic polynomial with {\em all} the roots \eqref{eq:roots}, and observe (e.g., by differentiating both sides and substituting $z=a+b$) that
\[
\Omega_n'(a+b)= \Pi_n^\text{ellipse}(a+b);
\]
thus the problem is transformed to the evaluation of $\Omega_n'(a+b)$. Then, construct the new polynomial $P_n(z) = \Omega_n(z) + (a^n+b^n)$, and note that it is uniquely characterized by the property that it is monic and equal to $a^n+b^n$ at all of the $z$-values \eqref{eq:roots}. 

Now, define yet a new polynomial $\mathcal{P}_n(z)$ via the initial data $\mathcal{P}_0(z) = 2$, $\mathcal{P}_1(z)=z$, and the recursive formula
\begin{equation}\label{eq:recursion}
\mathcal{P}_n(z) = z\mathcal{P}_{n-1}(z) - ab\mathcal{P}_{n-2}(z)
\end{equation}
for $n\geq 2$. It is immediate by induction on $n$ that the $\mathcal{P}_n(z)$ so constructed is monic of degree $n$. Then, prove by induction on $n$ that, for any $\theta\in [0,2\pi)$, this polynomial satisfies
\begin{equation}\label{eq:itheta}
\mathcal{P}_n(e^{i\theta}a + e^{-i\theta}b) = e^{in\theta}a^n + e^{-in\theta}b^n.
\end{equation}
This is verified for $n=0$ and $n=1$ by direct substitution of $e^{i\theta}a+e^{-i\theta}b$ for $z$, and then for higher $n$ by making this substitution in \eqref{eq:recursion}, invoking the induction hypothesis, and calculating:
\begin{align}
\mathcal{P}_n(e^{i\theta}a+e^{-i\theta}b) &= (e^{i\theta}a+e^{-i\theta}b)(e^{i(n-1)\theta}a^{n-1}+e^{-i(n-1)\theta}b^{n-1})\nonumber\\
&\quad\quad\quad-ab(e^{i(n-2)\theta}a^{n-2}+e^{-i(n-2)\theta}b^{n-2})\nonumber\\
&= e^{in\theta}a^n + e^{-in\theta}b^n.\label{eq:maincalc}
\end{align}

This proves \eqref{eq:itheta}, from which it immediately follows that, for $z$ of the form \eqref{eq:roots},
\[
\mathcal{P}_n(z) = a^n+b^n.
\]
This in turn means that $\mathcal{P}_n(z)$ obeys the properties that uniquely characterize $P_n(z)$, i.e., $P_n(z)=\mathcal{P}_n(z)$. Thus $P_n(z)$ is the polynomial defined by the recursion \eqref{eq:recursion} after all, and furthermore, it obeys \eqref{eq:itheta}. 

Finally, use \eqref{eq:itheta}, together with L'H\^{o}pital's rule, to carry out the evaluation of $\Pi_n^{\text{ellipse}}(a+b) = \Omega_n'(a+b)$:
\begin{align}
\Omega_n'(a+b) &= P_n'(a+b)\nonumber\\
&=\lim_{\theta \to 0} \frac{e^{in\theta}a^n+e^{-in\theta}b^n - (a^n+b^n)}{e^{i\theta}a+e^{-i\theta}b - (a+b)}\nonumber\\
&=\lim_{\theta \to 0} \frac{ine^{in\theta}a^n-ine^{in\theta}b^n}{ie^{i\theta}a-ie^{i\theta}b}\nonumber\\
&= n\frac{a^n-b^n}{a-b}.\label{eq:finalans}
\end{align}
This is the final answer. Price noted the similarity with the Binet formula for the Fibonacci numbers in \cite{price1}. In \cite{price2}, he noted that the Binet formula and \eqref{eq:finalans} together prove the pattern observed in Table \ref{tbl:ellipseproducts} after specializing to the appropriate values of $a,b$, and made the further observation that, in view of the initial data and the recursive formula \eqref{eq:recursion}, $P_n(z)$ is actually the generalized Lucas polynomial $V_n(X,Y)$ \cite{bergumhoggatt, swamy} with $X=z$ and $Y=-ab$.

We treat this resolution of the ellipse problem as the starting point for a further inquiry. The key devices of Price's argument---the shift of $\Omega_n(z)$ by $a^n+b^n$ to form $P_n(z)$, and its identification with the recursively defined polynomial $\mathcal{P}_n(z)$ via the calculation \eqref{eq:maincalc}---are effective but mysterious. Where are they coming from and why do they work? Furthermore, the fact that \eqref{eq:finalans} looks like the Fibonacci numbers' Binet formula, while $P_n(z)$ is a Lucas polynomial, suggests connections with known theory, but these connections are far from clear. 

Our primary goal in the rest of this article, carried out in Sections \ref{sec:cardano}--\ref{sec:solution}, is to give a reorganized version of the proof of \eqref{eq:finalans}, with a few new ideas, that reveals many points of contact with classical theory, providing motivation and illumination. Our method surfaces some additional links, not directly used in the proof of \eqref{eq:finalans}, and a secondary goal is to develop these. What follows is an outline of the rest of the article, giving the structure of the argument and the connections that are drawn.

As in Price's solution, the work lies in evaulating $\Omega_n'(a+b)$. We begin by showing in Section \ref{sec:cardano}, using an invariance argument, that the polynomial $\Omega_n(z)$ with the roots \eqref{eq:roots} is equal to $L_n(z,ab)-(a^n+b^n)$, where $L_n(X,Y)$ is a polynomial that expresses a relation between various standard symmetric polynomials in two variables. ($L_n(z,ab)$ coincides with Price's $P_n(z)$, though the definitions are different.) The argument is inspired by Cardano's solution of the cubic equation in radicals.

In Section \ref{sec:recursive}, we show that $L_n(X,Y)$ obeys a quadratic recurrence relation. We also give an explanation for this relation in terms of two broader mathematical contexts: the general theory of linear recurrences, and a classical theorem of Newton on power sums. 

In Section \ref{sec:binet}, we identify $L_n(X,Y)$ as the generalized Lucas polynomial $V_n(X,Y)$, up to a sign change. We give a parallel interpretation of the generalized Fibonacci polynomial. These are reformulations of the Binet formulas for these polynomials.

Section \ref{sec:solution} completes the evaluation of $\Omega_n'(a+b)$, and thus the determination of the product of the elliptical chord lengths, by applying a known relation between the generalized Lucas and Fibonacci polynomials. We also give a proof of this relation inspired by the L'H\^{o}pital's rule calculation in \eqref{eq:finalans}.

Relationships with Price's proof are discussed along the way.

Section \ref{sec:furtherremarks} develops an additional link surfaced by the argument of Sections \ref{sec:cardano}--\ref{sec:solution}. Since Section \ref{sec:cardano} reveals $\Omega_n(z)$ to be a generalization of Cardano's ``depressed cubic," we use Cardano's method to give a radical formula for the roots of $\Omega_n(z)$, which specializes to the known formula for the roots of Lucas polynomials.

\begin{remark}
The summary of Price's method for proving \eqref{eq:finalans} given in this section is not an exhaustive account of the work done in \cite{price1, price2}. In \cite{price1}, Price also considered the products of chord lengths that arise from rotating the roots of unity by a fixed angle along the unit circle prior to scaling. In \cite{price2}, he used the interpretation of Fibonacci and Lucas numbers in terms of products of elliptical chord lengths to recover identities among and divisibility properties of these numbers, such as $F_{2n} = F_nL_n$.
\end{remark}

\section{Symmetric polynomials and Cardano's depressed cubic.}\label{sec:cardano}

In this and the next three sections, we give our reformulation of the solution of the ellipse problem. Recall from the previous section that the task is to evaluate $\Omega_n'(a+b)$. The focus here is on carrying out the evaluation in such a way that for each step, a broader mathematical context is visible.

This section focuses on illuminating the polynomial $\Omega_n(z)$. We will see that as a consequence of an invariance property of its set of roots, $\Omega_n(z)$ has the form $L_n(z,ab) - (a^n+b^n)$, where $L_n(X,Y)$ is a polynomial bearing a natural definition in terms of the classical theory of symmetric polynomials.\footnote{Note that this means $L_n(z,ab)$ is equal to Price's polynomial $P_n(z)$. In spite of this equality, we are introducing the new notation here to emphasize that the definitions are {\em a priori} different.} We also discuss how the authors were led to this insight by the Cardano formula for the solution of a cubic equation in radicals.

The critical observation is that the set of roots \eqref{eq:roots} is the result of repeated applications of the substitution $a\mapsto \zeta a, b\mapsto \zeta^{-1}b$ to the expression $a+b$, and that this substitution leaves invariant the polynomials $ab$ and $a^n+b^n$. Therefore, any polynomial expression linking the three quantities $a+b$, $ab$, and $a^n+b^n$, will link $\zeta^ja+\zeta^{-j}b$, $ab$, and $a^n+b^n$ in the same way, for any $j$. In other words, any polynomial with coefficients in terms of $ab$ and $a^n+b^n$, and with $a+b$ as a root, will have {\em all} the complex numbers \eqref{eq:roots} as roots.

Such a polynomial relation between $a+b$, $ab$, and $a^n+b^n$ is guaranteed to exist by the fact that there are only two indeterminates; thus these three expressions are not algebraically independent. There will be many such polynomial relations; the one we want, $\Omega_n$, is monic and of minimal degree in $a+b$. We can bring it into clearer focus via the fundamental theorem on symmetric polynomials (FTSP). This classical theorem states that any polynomial in indeterminates $a_1,\dots, a_m$ that is symmetric, i.e., invariant under all permutations of the indeterminates, is representable in a unique way as a polynomial in the {\em elementary symmetric polynomials} $e_1 = \sum_j a_j$, $e_2 = \sum_{j<k} a_ja_k$, $e_3 = \sum_{j<k<\ell} a_{j}a_{k}a_{\ell}$, $\dots$. Furthermore, if the original symmetric polynomial to be represented has integer coefficients, then so does the polynomial representation in terms of the $e_i$'s.\footnote{This theorem was well known by the end of the 18th century. Its provenance is not straightforward to describe. H. Edwards argues that it was in essence known to Newton \cite[\S 8--9]{edwards1}, although he never published a proof---but see the special case discussed in the next section. J.-P. Tignol \cite[p.~99]{tignol} credits the first proof to Waring's 1770 {\em Meditationes Algebraicae}. The classic proof is found in an 1816 paper of Gauss \cite[paragraphs 3--5]{gauss2}. For an in-depth discussion of this proof, see \cite{BBSCoskey}. Gauss's paper also appears to contain the first modern {\em statement} of the theorem; all prior articulations formulated it as a statement about roots of a polynomial, implicitly assuming the existence of a splitting field. The theorem plays a foundational role in Galois's development of the theory of equations for which he is remembered; see \cite{edwards1} and \cite{tignol}. Proving it for two indeterminates---the case relevant to this article---is a nice exercise for the interested reader.}

In the present case of the two indeterminates $a,b$, the FTSP states that any polynomial unchanged by swapping $a$ and $b$ is representable uniquely as a polynomial in $a+b$ and $ab$. In particular, the {\em power sum} $a^n+b^n$ has such a representation. Thus we can make the following definition.

\begin{definition}\label{def:L}
For each $n\geq 0$, let $L_n(X,Y)$ denote the unique integer polynomial such that
\begin{equation}\label{eq:defofL}
L_n(a+b,ab) = a^n + b^n
\end{equation} is an identity. 
\end{definition}

We will see in Section \ref{sec:binet} that $L_n(X,Y)$ is actually the generalized Lucas polynomial up to a sign change, although the latter is not usually defined in this way. For now, we satisfy ourselves with the light this definition sheds on $\Omega_n(z)$:

\begin{prop}\label{prop:thepoly}
The polynomial
\[
L_n(z,ab) - (a^n + b^n)
\]
is the unique monic polynomial with the roots \eqref{eq:roots}, where $L_n$ is the polynomial of Definition \ref{def:L}.

In other words, $\Omega_n(z) = L_n(z,ab) - (a^n+b^n)$.
\end{prop}

\begin{proof}Since $a^n+b^n$ is homogeneous of degree $n$ in the indeterminates $a,b$, the expression $L_n(a+b,ab)$ must also be homogeneous of degree $n$. Thus $L_n(X,Y)$ is homogeneous if $X$ and $Y$ are given the weights $1$ and $2$, respectively. Furthermore, $X^n$ is the only monomial in $X,Y$ that, upon substituting $X\mapsto a+b, Y\mapsto ab$ and expanding, produces a polynomial containing $a^n$ or $b^n$, and it contains them both with coefficient $1$. Thus $X^n$ occurs in $L_n(X,Y)$ with coefficient $1$. All the other monomials of $L_n(X,Y)$, being of total degree $n$, are of degree less than $n$ in $X$. It follows that $L_n(z,ab) - (a^n+b^n)$ is monic of degree $n$ in $z$.

Now the polynomial $L_n(z,ab) - (a^n+b^n)$ has $a+b$ as a root by the defining property \eqref{eq:defofL} of $L_n(X,Y)$. Furthermore, since $ab$ and $a^n+b^n$ are both invariant under the substitution 
\begin{equation}\label{eq:sub}
a\mapsto \zeta a, b\mapsto \zeta^{-1}b,
\end{equation}
the polynomial $L_n(z,ab)-(a^n+b^n)$ is as well, and by repeated application of \eqref{eq:sub} to the root $a+b$, we conclude that it has all $n$ of the numbers \eqref{eq:roots} as roots. Since it is of degree $n$, we conclude these are all its roots. Thus, it fulfills the properties that uniquely characterize $\Omega_n(z)$, and we conclude $\Omega_n(z) = L_n(z,ab) - (a^n+b^n)$, as desired.
\end{proof}

The authors' inspiration for Proposition \ref{prop:thepoly} was the observation that, in the $n=3$ case, the polynomial whose roots have the shape \eqref{eq:roots} is precisely the {\em depressed cubic equation}
\begin{equation}\label{eq:cubic}
z^3 + pz + q = 0.
\end{equation}
In this context, the word ``depressed" refers to the fact that \eqref{eq:cubic} is monic with no quadratic term. The classical Cardano solution of the cubic polynomial in radicals begins by using a change of variables to get the general cubic into this form. Once this is achieved, the roots of the depressed cubic are given as the numbers \eqref{eq:roots}, where
\begin{equation}\label{eq:aandb}
a,b = \sqrt[3]{-\frac{q}{2}\pm\sqrt{\frac{q^2}{4}+\frac{p^3}{27}}}.
\end{equation}
(The cube roots must be chosen so that $ab = -p/3$.) Thus, the polynomial $\Omega_n(z)$ is a generalization of Cardano's depressed cubic.

To understand why the roots of the depressed cubic have the right shape, the authors considered the following interpretation of Cardano's formula, given by Ronald Solomon in \cite[pp.~50--51]{solomon}. The identity
\begin{equation}\label{eq:identity}
(a+b)^3 - 3ab(a+b) - (a^3 + b^3) = 0
\end{equation}
becomes the depressed cubic \eqref{eq:cubic} if $z=a+b$, $p = -3ab$, and $q=-(a^3 + b^3)$. Thus, finding $z$ satisfying \eqref{eq:cubic} can be achieved by finding $a$ and $b$ satisfying $ab = -p/3$ and $a^3 + b^3 = -q$, whereupon $z=a+b$ is the desired solution. This in turn can be accomplished by noting that $ab = -p/3$ implies $a^3b^3 = -p^3/27$; thus $a^3,b^3$ are quantities whose sum and product are known. It follows that $a^3,b^3$ solve the known quadratic equation $X^2 + qX - p^3/27=0$, the {\em resolvent quadratic} of \eqref{eq:cubic}. Applying the quadratic formula to this equation and taking cube roots yields \eqref{eq:aandb}.

The insight relevant to the situation at hand is that we have a choice of cube roots for $a,b$ but must preserve the known relation $ab = -p/3$. Thus we can twist $a$ by a factor of $\zeta$, but not without twisting $b$ by $\zeta^{-1}$. This is why the ``original root" $a+b$ leads to all the roots \eqref{eq:roots}. The authors arrived at Proposition \ref{prop:thepoly} by seeking to generalize this situation to $n>3$. The polynomial $L_n(X,Y)$ is defined to capture the identity \eqref{eq:defofL} that generalizes \eqref{eq:identity}. (In the $n=3$ case, \eqref{eq:identity} yields that $L_3(X,Y)=X^3-3XY$.)

We close this section with two remarks about how the work done here relates to Price's solution given in Section \ref{sec:outline}.

First, while Price's solution does not explicitly mention invariance under a substitution, our proof of Proposition \ref{prop:thepoly} (and thus our whole approach) hinges on it. One may wonder whether the invariance of $ab$ and of $a^n+b^n$ under \eqref{eq:sub} are hiding somewhere in Price's argument. In fact, they both are---though in two different places. The invariance of $a^n+b^n$ under \eqref{eq:sub} is equivalent to Price's assertion, having proven that the recursively defined polynomial $\mathcal{P}_n(z)$ satisfies \eqref{eq:itheta}, that this means it is equal to $a^n+b^n$ at the $z$-values \eqref{eq:roots}. This conclusion, in turn, is used to identify $\mathcal{P}_n(z)$ with $P_n(z)$. Meanwhile, the invariance of $ab$ under \eqref{eq:sub} is a special case of the fact that $ab = (e^{i\theta}a)(e^{-i\theta}b)$ for any $\theta$, and this latter invariance property is behind all the cancellation in Price's key calculation \eqref{eq:maincalc}.\footnote{Incidentally, although the identity \eqref{eq:itheta} is central to Price's argument, it plays no direct role in our approach. However, we can recover it from the above by noting that Proposition \ref{prop:thepoly} implies that $L_n(a+b,ab)=P_n(a+b)$, and then substituting $a\mapsto e^{i\theta}a, b\mapsto e^{-i\theta}b$ in both this and in \eqref{eq:defofL}.}

Second, the reason for the utility of shifting $\Omega_n(z)$ upward by $a^n+b^n$, arriving at $P_n(z)$, which satisfies a nice recursion---remains a mystery in \cite{price1,price2}. We view Proposition \ref{prop:thepoly} as the beginning of a resolution: $\Omega_n(z)$ was something natural minus $a^n+b^n$ in the first place. In the next section we will complete this picture by showing why this natural thing should be expected to satisfy a recursion.

\section{A recursive formula.}\label{sec:recursive}

In this and the next section, we focus on illuminating $L_n(X,Y)$. Here we observe that $L_n(X,Y)$ obeys a recurrence relation, and put this relation in two broader contexts: the general theory of linear recurrences, and a classical theorem of Newton on power sums.

Since $L_n(X,Y)$ represents the power sum $a^n+b^n$ (in terms of $a+b$ and $ab$), a recursive formula for $L_n$ needs to capture the relationship between $a^n+b^n$ and lower-degree power sums. The identity
\begin{equation}\label{eq:newton2}
(a^n + b^n) = (a+b)(a^{n-1}+b^{n-1}) - ab(a^{n-2}+b^{n-2}),
\end{equation}
valid for all $n\geq 2$, articulates such a relationship. Making the substitution $a^j + b^j\mapsto L_j(a+b,ab)$ (for $j=n, n-1, n-2$) followed by $a+b \mapsto X$ and $ab \mapsto Y$, we obtain that

\begin{prop}\label{prop:lucasrec}
The polynomial $L_n(X,Y)$ satisfies the recursive formula
\[
L_n(X,Y) = X L_{n-1}(X,Y) - Y L_{n-2}(X,Y)
\]
for all $n\geq 2$. \qed
\end{prop}

This linear recurrence characterizes $L_n(X,Y)$ for all $n$ once $L_0(X,Y)$ and $L_1(X,Y)$ are known. Since $L_0(a+b,ab) = a^0 + b^0 =2$, we have $L_0(X,Y) = 2$, and since $L_1(a+b,ab) = a^1 + b^1 = a+b$, we have $L_1(X,Y) = X$.

The heart of Proposition \ref{prop:lucasrec} is the identity \eqref{eq:newton2}. Although one can view it as an isolated piece of high school algebra, our interest is in broader mathematical contexts into which it fits. We mention two.

The first is the general theory of linear recurrences. Let $(g_n)_{n\geq 0}$ denote a sequence. A {\em homogeneous linear recurrence of degree $m$ with constant coefficients} is a condition of the form
\[
\sum_{j=0}^m c_jg_{n-j} = 0,\;\text{for all }n\geq m,
\]
where $m$ is a natural number and the $c_j$ are constants with $c_0$ nonzero. The family of sequences that satisfy this condition is closed under term-by-term addition and multiplication by scalars, so it forms a vector space with respect to these operations. It is of dimension $m$ since $g_0,\dots,g_{m-1}$ can be freely chosen, whereupon the recursion determines $g_n$ for $n\geq m$. 
A geometric sequence $(a^n)_{n\geq 0}$ belongs to this family if and only if $a$ is a root of the polynomial $\sum_{j=0}^m c_nX^{m-j}$, the {\em characteristic polynomial} of the recurrence. (For a more general but still very efficient account of this theory, geared toward its connection with generating functions, see \cite[Section 4.1]{stanley1}.) 

This theory motivates the following approach to arriving at \eqref{eq:newton2}. The geometric sequences $(a^n)$ and $(b^n)$ both satisfy the recurrence whose characteristic polynomial is the quadratic with $a,b$ as roots; this is $X^2 - (a+b)X + abX$. It follows that any linear combination of $(a^n)$ and $(b^n)$, and in particular the sequence $(a^n+b^n)$, also satisfies this recurrence. This yields \eqref{eq:newton2}, but also the generalization
\begin{equation}\label{eq:2linrec}
xa^n + yb^n - (a+b)(xa^{n-1} + yb^{n-1}) + ab(xa^{n-2} + yb^{n-2}) = 0,\; n\geq 2,
\end{equation}
where $x,y$ are arbitrary. It is straightforward to further generalize to $m>2$ indeterminates $a_1,\dots,a_m$; the characteristic polynomial of the recursion is then the generic monic polynomial with $m$ roots, whose coefficients are, up to sign, the elementary symmetric polynomials $e_1,\dots,e_m$. We get
\begin{equation}\label{eq:linrec}
\sum x_ja_j^n - e_1\sum x_ja_j^{n-1} + e_2\sum x_ja_j^{n-2} -\dots \pm e_m\sum x_ja_j^{n-m} = 0,
\end{equation}
where the $x_j$'s are arbitrary.

A second broader context, in fact the one in which the authors viewed \eqref{eq:newton2} when developing the present argument, is recommended by the fact that $L_n(X,Y)$'s definition is as a representation of power sums in terms of elementary symmetric polynomials. There is a completely classical theorem that supplies an inductive method for generating such representations. This special case of the FTSP appeared in the written record before the general case. {\em Newton's theorem}, published in Newton's {\em Arithmetica Universalis} in 1707, states that, for any $m$ indeterminates $a_1,\dots,a_m$, the {\em power sums} $p_1 = \sum_j a_j$, $p_2 = \sum_j a_j^2$, $\dots$, and the elementary symmetric polynomials $e_1,\dots,e_m$, obey the following relation:
\[
p_n - e_1p_{n-1} + e_2p_{n-2} - \dots \pm e_{n-1}p_1\mp ne_n = 0.
\]
(See \cite[\S 8,12]{edwards1}.) By convention, $e_j=0$ for all $j>m$, so when $n\geq m$, this reduces to the specialization of \eqref{eq:linrec} in which all $x_j$'s are $1$. However, it also holds when $n< m$. Indeed, it holds in the $m\to \infty$ limit, the {\em ring of symmetric functions}.\footnote{The standard reference on this combinatorially and representation-theoretically important ring is \cite[Chapter I]{macdonald}. More recent treatments include \cite[Chapter 7]{stanley} and \cite[Chapter 4]{sagan}. For a very quick introduction, see \cite[Chapter 6]{fulton}.}

The identity \eqref{eq:newton2} is, of course, the case $m=2$. 

The work in this section sheds light on the main calculation \eqref{eq:maincalc} in Price's proof. We view this and the previous section as providing a ``direct sum decomposition" of this calculation into (i) the invariance of $ab$ under $a\mapsto e^{i\theta}a, b\mapsto e^{-i\theta}b$, as discussed in Section \ref{sec:cardano}, and (ii) the identity \eqref{eq:newton2}, in turn an output either of the theory of linear recurrences or of Newton's theorem, as discussed here.

\section{The Binet formulas.}\label{sec:binet}

Readers familiar with the generalized Fibonacci and Lucas polynomials (see \cite{bergumhoggatt, hoggattbicknell, swamy, webbparberry}) will recognize in Proposition \ref{prop:lucasrec} the recurrence defining these polynomials, up to the sign change $Y\mapsto -Y$. The first two polynomials $L_0 = 2$ and $L_1 = X$ coincide with the generalized Lucas polynomials. It is then immediate, since $L_0$ and $L_1$ do not involve $Y$, that the following holds.

\begin{prop}\label{prop:lucas}
The polynomial $L_n(X,Y)$ of Definition \ref{def:L} is the image of the $n$th generalized Lucas polynomial $V_n(X,Y)$ under the substitution $Y\mapsto -Y$.\qed
\end{prop}

Since in Section \ref{sec:cardano} we saw that $L_n(z,ab)$ coincides with Price's $P_n(z)$, Propositions \ref{prop:thepoly}--\ref{prop:lucas} provide some illumination of Price's observation that $P_n(z)$ is a generalized Lucas polynomial. An explicit formula for the coefficients of $V_n(X,Y)$ is known, e.g., \cite[equation (2.22)]{swamy}.

Proposition \ref{prop:lucas} interprets the generalized Lucas polynomial in terms of symmetric polynomials. It is actually a reformulation of a standard fact about the generalized Lucas polynomials, their so-called {\em Binet formula}. In this section, we discuss this connection, give a parallel interpretation of generalized Fibonacci polynomials, based on an identity parallel to \eqref{eq:newton2}, and contextualize this identity in parallel with what was done for \eqref{eq:newton2} above.

The Binet formula for generalized Lucas polynomials (e.g., \cite[equation (14)]{bergumhoggatt} or \cite[equation (2.2)]{swamy}) asserts that
\begin{equation}\label{eq:binet}
V_n(X,Y) = a^n + b^n,
\end{equation}
where $a=(X + \sqrt{X^2 + 4Y})/2$ and $b=(X-\sqrt{X^2+4Y})/2$, or, equivalently, $X = a+b$ and $Y = -ab$. Normally, one sees $X,Y, V_n$ as conceptually prior to $a,b$. The usual proof is either by a mechanical induction that comes down to the identity \eqref{eq:newton2}, or via a standard exercise in finding a closed form for the linear recurrence $V_n = XV_{n-1} + YV_{n-2}$. In the latter case, $a,b$ enter as the roots of the characteristic polynomial (see the previous section) of that recurrence. However, in our present context, which treats $a,b$ as the conceptual starting point (as in Definition \ref{def:L}), the Binet formula \eqref{eq:binet} asserts that $V_n(X,Y)$ is the polynomial that expresses the power sum $a^n+b^n$ in terms of the (sign-adjusted) elementary symmetric polynomials $X = a+b$ and $Y = - ab$. In view of Definition \ref{def:L}, this is the content of Proposition \ref{prop:lucas}.

We can give a parallel interpretation for the generalized Fibonacci polynomials. The latter, denoted $U_n(X,Y)$, are defined by the same recurrence as the Lucas polynomials, but with the initial data $U_0 = 0$, $U_1= 1$. 

The Binet formula for generalized Fibonacci polynomials (\cite[Theorem 1]{webbparberry} or \cite[equation (2)]{hoggattbicknell}) states that
\begin{equation}\label{eq:fibobinet}
U_n(a+b,-ab) = \frac{a^n-b^n}{a-b}.
\end{equation}
The right side is equal to 
\[
h_{n-1} = a^{n-1}+a^{n-2}b+\dots+ab^{n-2} + b^{n-1},
\]
the {\em complete symmetric polynomial} of degree $n-1$ in the two indeterminates $a,b$. This is defined to be the sum of all monomials of the given degree. Like the power sums $p_j$ and the elementary symmetric polynomials $e_j$, the complete symmetric polynomials $h_j$ are basic objects in the theory of symmetric functions.\footnote{In fact, they are algebraically interchangeable with the elementary symmetric polynomials: there is an involution of the ring of symmetric functions that exchanges these two classes (\cite[Theorem 7.6.1]{stanley} or \cite[Section 6.2, Corollary 1]{fulton}). This is a consequence of the FTSP and the identity \eqref{eq:completehom} below.} Thus, the following definition is a natural one.

\begin{definition}\label{def:F}
Let $F_n(X,Y)$ be the integer polynomial such that
\[
F_n(e_1,e_2) = h_{n-1}
\]
is an identity, where $e_1, e_2$ are the elementary symmetric polynomials in two indeterminates $a,b$, and $h_{n-1}$ is the complete symmetric polynomial of degree $n-1$.
\end{definition}

With this definition, $F_n$ is to $U_n$ as $L_n$ is to $V_n$. With $X=e_1=a+b$ and $Y=e_2=ab$, we have
\begin{align*}
F_n(X,Y) &= h_{n-1}\\
&= U_n(a+b,-ab)\\
&= U_n(X,-Y).
\end{align*}
The first equality is by definition of $F_n$, and the second is the Binet formula. So Definition \ref{def:F} interprets the generalized Fibonacci polynomial in terms of symmetric polynomials in the same way Definition \ref{def:L} does this for the generalized Lucas polynomial.

The proof of Proposition \ref{prop:lucas}, or, equivalently, the Binet formula for generalized Lucas polynomials, comes down to the identity \eqref{eq:newton2}. In a similar way, the standard inductive proof of the Binet formula for generalized Fibonacci polynomials comes down to the identity
\begin{equation}\label{eq:complete2}
\frac{a^n-b^n}{a-b} - (a+b) \frac{a^{n-1}-b^{n-1}}{a-b} + (ab)\frac{a^{n-2}-b^{n-2}}{a-b} = 0.
\end{equation}
(Alternatively, in a more precise parallel to what was done here for $V_n(X,Y)$, one can derive the equality $F_n(X,Y)=U_n(X,-Y)$ directly from this identity, and then reinterpret this latter equality {\em as} the Binet formula.) 

In the previous section, we contextualized \eqref{eq:newton2} both as an output of the theory of linear recurrences and as the $m=2$ case of Newton's theorem. We now seek to similarly contextualize \eqref{eq:complete2}. It comes from the theory of linear recurrences in exactly the way \eqref{eq:newton2} did; indeed, it is just \eqref{eq:2linrec} in the special case $x =(a-b)^{-1}, y=-(a-b)^{-1}$. But is there an analogue to Newton's theorem?

There is. For any number of indeterminates $m$, the $e_j$'s and $h_j$'s satisfy a relation \cite[Ch. I, equation (2.6')]{macdonald} nearly identical to the one between the $e_j$'s and the $p_j$'s captured by Newton, namely
\begin{equation}\label{eq:completehom}
h_n-h_{n-1}e_1+h_{n-2}e_2 - \dots \pm he_{n-1} \mp e_n = 0
\end{equation}
for all $n\geq 1$. As in Section \ref{sec:recursive}, this continues to hold for $n>m$ via the convention that $e_j=0$ for $j>m$. The standard proof is a beautiful and very short calculation with generating functions; see \cite[p.~21]{macdonald}, \cite[p.~296]{stanley}, or \cite[p.~73]{fulton}. In analogy with \eqref{eq:newton2} with respect to Newton's theorem, \eqref{eq:complete2} is the $m=2$ case.

\section{Stretched chords.}\label{sec:solution}

In this section we apply the insights about $\Omega_n(z)$ and $L_n(X,Y)$ obtained in Sections \ref{sec:cardano}--\ref{sec:binet} to evaluate $\Omega_n'(a+b)$, completing our reformulation of the solution to the ellipse problem. The key device is a known relation \cite[equation (3.10)]{swamy} between generalized Lucas and Fibonacci polynomials:
\begin{equation}\label{eq:UandV}
\frac{d}{dX} V_n(X,Y) = n U_n(X,Y).
\end{equation}
We also give a proof of this relation inspired by the L'H\^{o}pital's rule calculation in Price's argument.

The previous two sections yield that
\[
\Omega_n'(z) = \frac{d}{dz} L_n(z,ab) = \frac{d}{dz} V_n(z,-ab).
\]
(The first equality follows from Proposition \ref{prop:thepoly}, and the second from Proposition \ref{prop:lucas}.) Applying \eqref{eq:UandV} (with $X=z, Y=-ab$) to the right side, and then evaluating at $z=a+b$ on both sides, we get
\[
\Omega_n'(a+b) = nU_n(a+b,-ab).
\]
Combining this with the Binet formula for generalized Fibonacci polynomials \eqref{eq:fibobinet}, we recover Price's beautiful result:

\begin{prop}[Price \cite{price1}]\label{prop:price}
The product of the elliptical chord lengths described in the introduction is
\[
\pushQED{\qed}
n\frac{a^n-b^n}{a-b}.\qedhere
\popQED
\]
\end{prop}

The identity \eqref{eq:UandV} used in this section plays a parallel role to that played in Price's argument by the L'H\^{o}pital's rule calculation in \eqref{eq:finalans}---both allow passage from $\Omega_n(z)$ to its derivative at $z=a+b$. In \cite{swamy}, the relation \eqref{eq:UandV} is proven indirectly via generating functions. However, this parallel suggests the possibility that \eqref{eq:UandV} also has a direct proof similar to the L'H\^{o}pital calculation in \eqref{eq:finalans}. Here is one such proof.

Set $a+b=X$ and $ab=-Y$ as usual. Treating $Y$ as a constant, we find that $b$ and therefore $X$ and $V_n(X,Y)$ become rational functions of $a$. We have $b=-Ya^{-1}$, so that $db/da = Ya^{-2} = -ba^{-1}$. Thus, by the Binet formulas \eqref{eq:binet} and \eqref{eq:fibobinet} for $V_n$ and $U_n$ and the chain rule, we have
\begin{align*}
\frac{d}{dX} V_n(X,Y) &= \frac{dV_n(X,Y)/da}{dX/da}\\
&= \frac{d(a^n + b^n)/da}{d(a+b)/da}\\
&= \frac{na^{n-1} +nb^{n-1}(-ba^{-1})}{1-ba^{-1}}\\
&= n\frac{a^n - b^n}{a-b}\\
\pushQED{\qed}
&= nU_n(X,Y).\qedhere
\popQED
\end{align*}

Compare with \eqref{eq:finalans}.

\section{Solution of $\Omega$ by radicals.}\label{sec:furtherremarks}

The exposition in Sections \ref{sec:cardano}--\ref{sec:solution} aimed to situate the ellipse problem in a dense web of connections. Now that the solution is complete, we conclude this article by developing one final connection. 

It was noted in Section \ref{sec:cardano} that the key polynomial $\Omega_n(z)$ generalizes Cardano's depressed cubic. We take this as an invitation, and solve $\Omega_n(z)$ in radicals using Cardano's method. We recover the known formula for the roots of the Lucas polynomials as a special case.

For $n=3$, we have
\[
\Omega_3(z) = z^3 - 3abz - (a^3 +b^3)
\]
by comparing \eqref{eq:identity} with Definition \ref{def:L} and Proposition \ref{prop:thepoly}. In Section \ref{sec:cardano}, we identified this with Cardano's depressed cubic $x^3 + px + q$ by setting $-3ab = p$ and $-(a^3+b^3)=q$. The presumption is that $p$ and $q$ belong to a prespecified field, such as the rational numbers. Rationality would not have been affected by instead setting $ab = p$ and $a^3 + b^3=q$, so moving forward, this is the convention we generalize: let $p=ab$ and let $q=a^n+b^n$. Then $\Omega_n(z)=L_n(z,p)-q$ is a univariate polynomial over a field that contains $p$ and $q$, and the goal of a solution in radicals is an expression for its roots in terms of $p$ and $q$.

We already know the roots are given by \eqref{eq:roots}; thus we only need to express $a,b$. Following Cardano (see Section \ref{sec:cardano}), we know that $a^nb^n=p^n$, thus, $a^n$ and $b^n$ are roots of the quadratic equation
\[
X^2 - qX + p^n = 0.
\]
Applying the quadratic formula and taking $n$th roots, we obtain
\begin{equation}\label{eq:solnbyradicals}
a,b = \sqrt[n]{\frac{q \pm \sqrt{q^2-4p^n}}{2}}, 
\end{equation}
where the choice of roots must respect the equality $ab=p$. So a radical expression for the roots of $\Omega_n(z) = L_n(z,p)-q$ is obtained by substituting \eqref{eq:solnbyradicals} in \eqref{eq:roots}.

The generalized Lucas polynomials are $V_n(X,Y) = L_n(X,-Y)$, i.e., the case $p=-Y$, $q=0$. In this case, \eqref{eq:solnbyradicals} simplifies to 
\[
a,b = \sqrt[n]{\pm\sqrt{-(-Y)^n}} = \sqrt[n]{\pm i\sqrt{-Y}^n} = \xi\sqrt{-Y}, \xi^{-1}\sqrt{-Y},
\]
where $\xi$ is an $n$th root of $i$, i.e., a $4n$th root of unity that is not a $2n$th root. We can take $\xi = e^{\pi i / 2n}$, $\zeta = \xi^4$, for example, and then by \eqref{eq:roots}, the roots of $V_n(X,Y)$ (as a polynomial in $X$) are
\[
\sqrt{-Y} (\zeta^j\xi + \zeta^{-j}\xi^{-1})= 2\sqrt{-Y} \cos \frac{(1+4j)\pi}{2n}
\]
for $j=0,1,\dots,n-1$. The expression $(1+4j)\pi/2n$ ranges over all the $k$-multiples of $\pi / 2n = 2\pi / 4n$, where $k=1\operatorname{mod} 4$, in the range $[0,2\pi)$. Exploiting the symmetry $\cos \theta = \cos(2\pi - \theta)$, those in the range $[\pi,2\pi)$ can be replaced by the $\ell$-multiples, where $\ell = 3\operatorname{mod} 4$, in the range $[0,\pi)$. Thus a slight simplification is
\[
X = 2\sqrt{-Y}\cos\frac{(1+2j)\pi}{2n}
\]
for $j=0,1,\dots,n-1$. Compare \cite[equation (2.24)]{swamy}. The classical Lucas polynomials are the case $Y=1$, so $\sqrt{-Y} = i$. Compare \cite[p.~274]{hoggattbicknell}.

\begin{acknowledgment}{Acknowledgments.}
Several individuals were involved in calling the authors' attention to Price's work: Francis Su, who included the special case of Table \ref{tbl:ellipseproducts} in his Harvey Mudd Math Fun Facts \cite{su}; Bowen Kerins and Darryl Yong, who then included the problem of computing this product in a problem set during a summer course at the Park City Mathematics Institute; and Sam Shah, who posted the problem on his blog \cite{shah}, whereby the authors became acquainted with it. The authors wish to thank Kerins in particular, who was a generous correspondent, and also Thomas Price himself, Tom Edgar, and Harold Edwards, for useful comments and for alerting us to \cite{price3}, \cite{galovich}, and \cite{kummerw}, respectively. We would also like to thank the anonymous referees, whose thoughtful feedback greatly improved the article.
\end{acknowledgment}

\begin{biog}
\item[Ben Blum-Smith] received his Ph.D. in 2017 from the Courant Institute of Mathematical Sciences at NYU, after a decade as a teacher and teacher educator. He is currently a Visiting Academic at the NYU Center for Data Science, and has also taught courses at Eugene Lang College, the Bard Prison Initiative, and the Bard Master of Arts in Teaching program. His research interests lie in invariant theory, algebraic combinatorics, their applications to data science, and connections between mathematics and democracy.
\begin{affil}
NYU Center for Data Science, New York, NY 10011\\
ben@cims.nyu.edu
\end{affil}

\item[Japheth Wood] received his Ph.D. in 1997 from the University of California, Berkeley. He has had academic appointments at the Pontificia Universidad Cat\'{o}lica de Chile, Vanderbilt University, the University of Louisville and Chatham University. Since 2006 he has been a faculty member at Bard College in Annandale-on-Hudson, New York, including stints with the Master of Arts in Teaching program and the Bard Prison Initiative. Japheth led the New York Math Circle from 2012 to 2015 and has taught at the Hampshire College Summer Studies in Mathematics program. He is co-director of the Bard Math Circle.
\begin{affil}
Mathematics, Bard College, Annandale-on-Hudson, NY 12504\\
jwood@bard.edu
\end{affil}
\end{biog}
\vfill\eject

\end{document}